# Fortified Test Functions for Global Optimization and the Power of Multiple Runs

Charles Jekel and Raphael T. Haftka

Department of Mechanical and Aerospace Engineering, University of Florida

**Abstract:** Some popular functions used to test global optimization algorithms have multiple local optima, all with the same value. That is all local optima are also global optima. This paper suggests that such functions are easily fortified by adding a localized bump at the location of one of the optima, making the functions more difficult to optimize due to the multiple competing local optima. This process is illustrated here for the Branin-Hoo function, which has three global optima. We use the popular Python SciPy differential evolution (DE) optimizer for the illustration. DE also allows the use of the gradient-based BFGS local optimizer for final convergence. By making a large number of replicate runs we establish the probability of reaching a global optimum with the original and fortified Branin-Hoo. With the original function we find 100% probability of success with a moderate number of function evaluations. With the fortified version, the probability of getting trapped in a non-global optimum could be made small only with a much larger number of function evaluations. However, since the probability of ending up at the global optimum is usually 1/3 or more, it may be beneficial to perform multiple inexpensive optimizations rather than one expensive optimization. Then the probability of one of them hitting the global optimum can be made high. We found that for the most challenging global optimum, multiple runs reduced substantially the extra cost for the fortified function compared to the original Branin-Hoo.

I. **Introduction**

There is large variety of global optimization algorithms, including nature inspired stochastic algorithms like simulated annealing, genetic algorithms, differential evolution (DE), particle swarm optimization, and ant colony optimization. There are also deterministic global optimizers like DIRECT (Jones et al., 1993) and SHGO (Enders et al.,2018). Finally, there are surrogate based adaptive sampling algorithms such as EGO (Jones et al., 1998).

Many test functions are used for tuning these algorithms, and a substantial list may be found in Wikipedia https://en.wikipedia.org/wiki/Test_functions_for_optimization. Unfortunately, few test functions have local optima that are competitive with the global optimum. Therefore, the chance that the algorithm will be trapped in a local optimum is low. Some popular test functions, like the Branin-Hoo function, and the Himmelblau function present have multiple local optimal, but all with the same value. That is, all their local optimal are global! Hence there is no chance of being trapped in a local non-global optimum, and finding a global optimum is easy because there are several to choose from.

One objective of this paper is to show that it is possible to convert these functions to more challenging problems with the global optimum being accompanied by competing local optima. This is done by adding a radial-basis-function (RBF) bump to one of the global optima. The Branin-Hoo function is used to illustrate how this makes the optimization much more challenging for a popular Python differential evolution optimizer in SciPy. A second objective is to explore how performing multiple short runs instead of one long one, can reduce the additional cost associated with the fortified function.

II. **The Branin-Hoo function**

Using the information from https://www.sfu.ca/~ssurjano/branin.html the Branin-Hoo function is defined as

$$f(\mathbf{x}) = a\left(x_2 - bx_1^2 + cx_1 - r\right)^2 + s(1-t)\cos(x_1) + s \ . \tag{1.1}$$

We use the recommended values $a=1, b=5.1/(4\pi^2), c=5/\pi, r=6, s=10, t=1/(8\pi)$.

The domain is here as usual, $x_1 \in [-5,10]$, $x_2 \in [0,15]$, and the function is shown in Fig. 1.

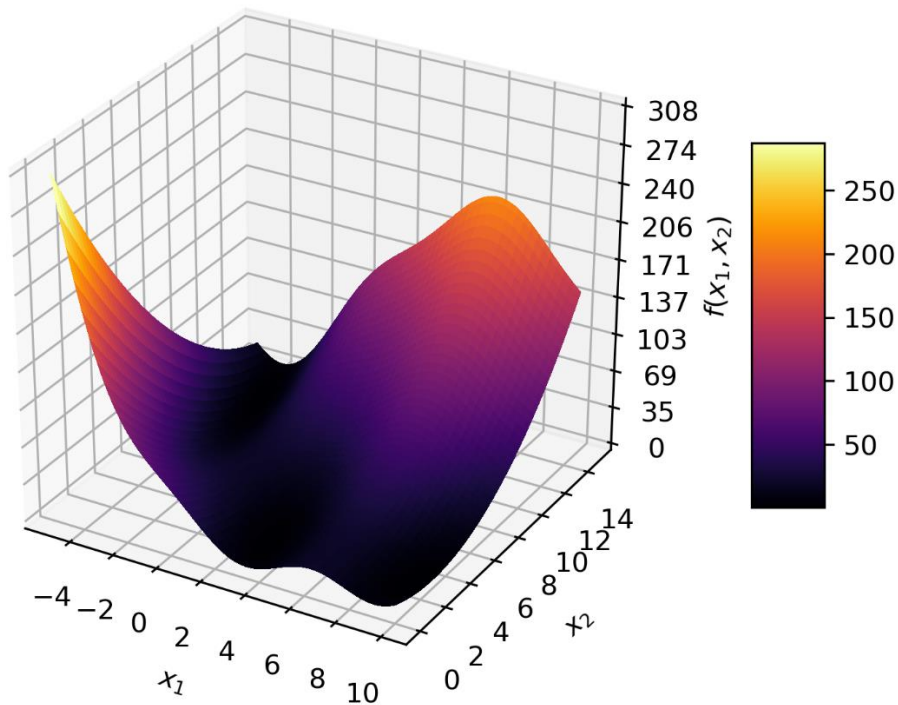

*Figure 1 Branin-Hoo function*

The function has three global optima where its value *is f(x)=0.397887*. They are numbered here as

First: $\mathbf{x} = (-\pi, 12.275)$

Second: $\mathbf{x} = (\pi, 2.275)$

Third: $\mathbf{x} = (9.42478, 2.475)$

III. **The bump**

It is desirable to add (for maximization) or subtract (for minimization) a bump to the original function that will not change the location of an optimum. Radial basis functions (RBF, Broomhead and Lowe,

1988) are selected because they depend only on the distance from the optimum. It is desirable that the bump will affect only one optimum, and for that the RBF bump function https://en.wikipedia.org/wiki/Radial_basis_function is selected. The bump function is defined as

$$\varphi(r) = \begin{cases} \exp\left(-\dfrac{1}{1-(\varepsilon r)^2}\right) & \text{for } r < \dfrac{1}{\varepsilon} \\ 0 & \text{otherwise} \end{cases} \quad (1.2)$$

Here *r* is the radial distance from the center of the bump, and its maximum value at *r=0* is *1/e*=0.3679. The width of the bump is determined by $\varepsilon$. Figure 2 shows a one-dimensional slice of the Branin-Hoo function with $10\varphi(r)$ subtracted at the location of its first global optimum.

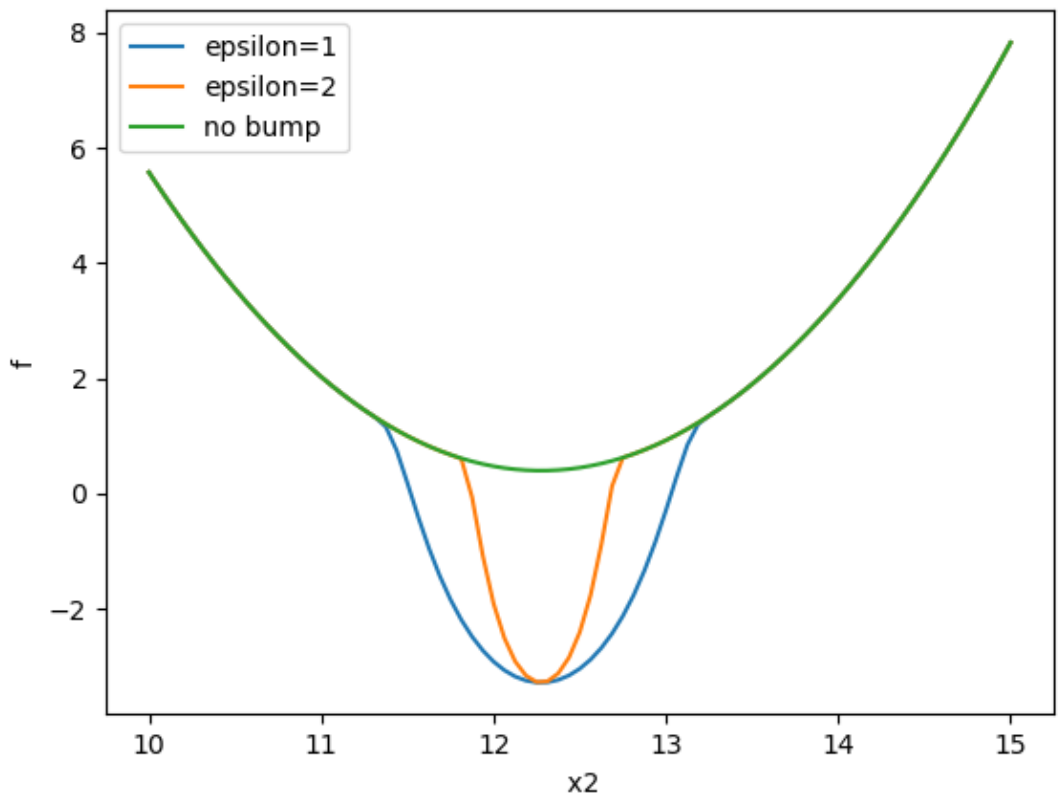

*Figure 2: A 1D slice through the Branin-Hoo function for $x_1 = -\pi$ with the bump defined by Eq. 1 multiplied by 10 subtracted at the location of its first global optimum at $(-\pi, 12.275)$. The width parameters are $\varepsilon$ = 1 and $\varepsilon$ = 2.*

IV. **Differential Evolution algorithm with follow-up by BFGS in SciPy**

Differential evolution (DE) (Storn and Price 1997) is a popular nature inspired global optimization algorithm. The algorithm is frequently used as a black-box optimizer, because DE does not require the objective function to be differentiable or continuous. The implementation of DE used here is available in the Python SciPy ecosystem (Virtanen et al., 2019). This implementation was chosen for rich features

(parallelization, Latin hypercube sampling initial design, L-BFGS-B follow up, etc.), simplicity, and wide availability.

Like most global search algorithms, once the algorithm gets close to the global optimum the final convergence is slow. The SciPy DE function therefore offers the option of switching to a local gradient based optimizer, BFGS, for this 'polish' phase. The number of function evaluations in the DE search is determined by the maximum number of iterations, *max_iter* and the population size *pop.* Once this phase is finished, the BFGS algorithm is run to convergence. Note that by specifying a small population size and a small maximum number of iteration the SciPy function will be using BFGS for most of the search.

**Performance measurement**

We define success for a run as reaching a global optimum with some given tolerance and failure as missing this target. DE is a stochastic algorithm, so that every time it is run one can get a different result. Additionally the initial design is randomized with a Latin hypercube sampling. Therefore, with the same algorithmic setting (e.g. number of iterations) some runs may fail and some succeed. When comparing the performance with and without a bump, we make multiple runs and compare the probability of failure estimated from the percentage of runs that fail Here we need to choose the number of runs, which is designated as *n* needed for desired accuracy in the probability estimate from the percentage of runs that fail.

The accuracy of the probability of failure estimate is measured here by the standard deviation of the number of runs that failed. It is easy to show that the standard deviation $\sigma_{fail}$ of the number of failures in *n* runs is related to the probability of failure *p* as

$$\sigma_{fail} = p\sqrt{pn(1-p)}. \tag{1.3}$$

To achieve accuracy of one percent, one sets $\sigma_{fail} = 0.01n$ and calculate the required number of runs for a given p.

$$(0.01n)^2 = p^3 n(1-p) \Rightarrow n = 10,000 p^3 (1-p) \tag{1.4}$$

While *p* is not known ahead of time, the maximum required *n* is attained for *p=0.75*, which gives *n=1055* runs. Since p=0.75 is the worst case, we rounded n to 1,000 for ease of translating the number of failures to percentages.

**Performance with no bump**

The DE algorithm was run first with the original Branin-Ho function for large number of times with for different population sizes (*pop*) and different numbers of maximum iterations (*max_iter*) with and without the BFGS polish. No convergence criteria were applied, so the optimization stopped with a fixed number of function evaluations equal to *2\*pop\*(max_iter+1)*. The optimization run was deemed successful if the optimized function value was higher than the true optimum by no more than 0.01. We recorded the percentage *p* of failed runs, as well as how many were at or near (within a Euclidean distance of *r=1*) from each optimum. This will provide insight on the behavior with a bump, discussed in

Section V. Note that the closest optima, Optimum 2 and Optimum 3 are at a distance of about 6.3 from each other.

We also made the same runs with the BFGS polish turned on to improve local convergence. This time all runs were successful, since all optima have the same value, but again we recorded how many were at each optimum. To make the runs repeatable, a random seed was specified.

Table 1 provides a summary of the results. As expected by Equation 1.4 and the discussion following the equation, with 1,000 runs the percentage results were mostly repeatable with a different random seed to about 1-2% accuracy. That is, if the percentage of failures is reported for example as 46%, when repeated with a different random seed it may range from 44% to 48%.

The first row in Table 1 shows that we can reach very low probability of missing a global optimum with 420 function evaluations, as only one optimization out of 1,000 failed. We also see that the DE favors the centrally positioned Optimum 2 compared to the other two. The second row shows that if we reduce the number of iterations to 10, almost half of the runs fail, but they all get within a distance of 1 from an optimum. That is, the failure is in the local convergence. Failed local convergence is supported by row 3, which shows that no failures occur if we add the BFGS follow-up. The average number of function evaluations is 240, with 220 coming from the DE part and 20 from the BFGS follow-up. With a population of 5 and maximum number of iterations equal to 5, 98% of the runs fail if we use only DE. Still 74% (15+37+22) make it to the vicinity of one of the optima. Adding BFGS reduces the number of failures to zero for an average cost of 84 function evaluations. Sixty of these come from the DE and 24 from the BFGS. This still may be viewed as a case where the two algorithms share the work.

The last row in Table 1 looks at a case where we use the BFGS to do practically the entire optimization. With a population of 2 and max number of iterations equal to 2 all runs fail, and only 23% get near one of the optima. When we add BFGS, they all converge with an average cost of 42 function evaluations. This time we see that when BFGS does most of the work the three optima have similar probabilities of being selected, while Optimum 2 was preferred in the previous runs where DE performed most of the work.

*Table 1 Performance of algorithms for the original Branin-Hoo function without a bump. Results are based on 1,000 replicate runs, and repeating the runs with different random seed indicates that the numbers are good to about 1% accuracy as predicted by Eqs. 1.3 and 1.4.*

| Algorithm | Pop | Max_iter | Percent failures | Average number of function evaluations | Percentage at or near each optimum |
|---|---|---|---|---|---|
| DE | 10 | 20 | 0.1% | 420 | 30, 44, 25 |
| DE | 10 | 10 | 46% | 220 | 27, 49, 24 |
| DE/BFGS | 10 | 10 | 0% | 240 | 27, 49, 24 |
| DE | 5 | 5 | 98% | 60 | 15, 37, 22 |
| DE/BFGS | 5 | 5 | 0% | 84 | 27, 45, 29 |
| DE | 2 | 2 | 100% | 15 | 7, 10, 6 |
| DE/BFGS | 2 | 2 | 0% | 42 | 30, 35, 35 |

V. **Performance with Bump at the first Optimum**

*The procedure described with no bump was repeated with the wider bump shown in Fig. 2, that is for an amplitude of* 10/e, *with* $\varepsilon = 1$. *The results when the bump was added to Optimum 1 are summarized in Table 2.Table 2 Performance of algorithms when a bump is added to Optimum 1 Results are based on 1,000 replicate runs and repeating the runs with different random seed indicates that the numbers are good to about 1% accuracy, as predicted by Eqs. 1.3 and 1.4.*

| Algorithm | Pop | Max iter | Percent failures | Average number of function evaluations | Percentage at or near each optimum |
|---|---|---|---|---|---|
| DE | 10 | 20 | 64 | 420 | 36, 39, 25 |
| DE | 20 | 10 | 89 | 440 | 46, 34, 20 |
| DE/BFGS | 20 | 10 | 54 | 458 | 46, 34, 20 |
| DE | 40 | 5 | 98 | 480 | 56, 27, 17 |
| DE/BFGS | 40 | 5 | 43 | 500 | 56, 27, 17 |
| DE | 50 | 4 | 99 | 500 | 57, 22, 21 |
| DE/BFGS | 50 | 4 | 43 | 522 | 57, 22, 21 |
| DE | 80 | 5 | 96 | 960 | 71, 14, 15 |
| DE/BFGS | 80 | 5 | 29 | 980 | 71, 14, 15 |
| DE/BFGS | 100 | 4 | 26 | 1021 | 74, 14, 12 |
| DE/BFGS | 125 | 3 | 19 | 1022 | 82, 8, 10 |
| DE/BFGS | 165 | 2 | 16 | 1013 | 84, 7, 9 |
| DE/BFGS | 330 | 2 | 4 | 2002 | 96, 2, 2 |

The first row of Table 1, with a population of 10 and 20 iterations with DE alone had only one failure out of a thousand runs when run without a bump. This was due to the fact that with 20 iterations DE can achieve final convergence, and it did not matter which optimum it converged to. With the bump, Table 2 shows 64% failures, because the percentage of runs that converged to Optimum 1 was 36%. Note that without the bump, as can be seen from Table 1, only 30% of the runs converged to Optimum 1.

Increasing the number of iterations, or adding BFGS does not help, because the runs are already converged, so we next increase the population size to 20 and reduced the number of iterations to 10 to get a similar cost. The second row in Table 2 shows that this increased the percentage of runs that went to the bump (Optimum 1) to 46%. The percent failure increased to 89% because 10 iterations were not close enough for final convergence. However, here adding BFGS (row 3) caused all of these 46% runs to converge to Optimum 1, reducing the percent failures to 54%.

Based on the results in Table 1, even 5 iterations were enough to get DE close to the optima, so we next tried a population of 40 with 5 iterations. This increased the percentage of runs near Optimum 1 to 56%. With DE alone 98% of the runs failed because of lack of final convergence, but adding BFGS we realized the 56% rate of success, or 44% failure. The average number of function evaluations, shows that adding BFGS increases the cost by about 20 evaluations, in line with what we saw without the bump in Table 1.

Additional experimentation showed that as the population size increases, the optimum number of iterations decreases. This is due to the two dimensional nature of the problem that creates a high probability that the initial population, placed randomly, will have a member inside the bump.

The total area of the design space is 15x15=225. The bump covers a circle of radius 1, so that its area is equal to $\pi$. The probability that a member of the population will be inside the bump area is $\pi/225$ or about 1.4%. The probability that it will outside is about 98.6%. However, with a population size of *pop*, the probability that not even one member of the population will fall inside the bump is $0.986^{pop}$. When *pop=165,* this probability is 9.8%, so there is more than 90% chance that one member of the initial population will be inside the bump, and several other members very close to the bump.

Comparing Table 1 to Table 2 demonstrates that adding the bump increases the required number of function evaluations by more than one order of magnitude.

VI.  **The power of multiple runs**

Differential evolution is a stochastic algorithm, which means that each run can produce different results. For such algorithms Schutte et al. (2007) suggested that it may be more efficient to run the algorithm multiple times, rather than let it run for a long time. We explore here the possibility of multiple runs with small population size instead of one run with large population size.

The concept is built on the assumption that the runs are at least approximately independent. That means that if the probability of failure in one run is denoted by *p*, the probability of failure in *m* multiple independent runs is

$$p_m^{ind} = p^m \tag{1.5}$$

This probability will be compared to the observed percentage $p_m$ of runs that fail.

Unfortunately, the number of runs required for good accuracy increases substantially. For example, if we select *m=3*, and we want 1,000 replications to estimate the probability of failure, we will need a total of 3,000 runs. If in addition we want to check whether the runs are approximately independent, we will need even more runs. This is because if the error in *p* is 1%, the error in $p^3$ is approximately 3%. So to achieve an error of 1% in $p^3$, we will need to reduce the error in p to about 1/3%. Following an analysis similar to Eq. 1.4, we find that we need approximately 9,000 runs to achieve the desired accuracy.

The procedure for checking the effectiveness of multiple runs and the accuracy of the independence check is first illustrated with small number of 30 runs. Since we plan to count the number of failures, we denote a successful run by "1" and a run that failed to converge close enough to the optimum by "0". To illustrate the difference between the true probability and the observed percentages, we chose p=2/3, and generated random integers with that probability. For one random seed (seed=30), we obtained 000001100000000000110110000101, counting 22 failures compared to the expected 20, so that the observed percentage of failures is 73.3% instead of 66.7%. For a different seed (seed=10) we obtained 000110100100010101011101000110 with 17 failures corresponding to 57% observed probability of failures. For p=2/3, and n=30, Eq. 1.3 gives a standard deviation of 1.72 failures, so that the observed variation is of the expected magnitude.

We now use these two cases to illustrate the effect of multiple runs by considering the case of m=3. That is, we assume that we will make three runs and take the best result of the three. In terms of success and failure, we will have failure only if all three runs failed. If we take the runs with seed=30 and we group them in triads, we get 000-001-100-000-000-000-110-110-000-101. Out of the 10 triads, five showed 000 corresponding to failure. Based on the true probability of failure of 2/3, we would have expected a percentage of failure of $100(2/3)^3 = 29.6\%$. Based on the observed 22 failures we would have expected $100(22/30)^3 = 39.4\%$. This is still substantially lower than the observed 50% but it illustrates the fact that errors magnify by about a factor of 3. Here a 10% difference between the observed number of failures (22) and the ideal one (20) magnified to about 33% difference when raised to the third power.

There is also the error due to the fact that we are estimating probability from only 10 groups of runs, so using Eq. 1.3 with p=0.296 (29.6%) gives a standard deviation of about 0.43 failures. For the seed of 10 we get 000-110-100-100-010-101-011-101-000-110 with only two failures out of 10.

Table 3 explores this option and the goodness of the independence assumption. For m=10, we would need about 100,000 runs to achieve the desired accuracy. We chose 100,800 runs because this number is divisible by all integers from 2 to 10. So, for example, if we specify that we check the probability of failure with a multiple of 7 runs, we divide the 100,800 to 14,400 groups of 7 runs each. This means that we can use these 100,800 runs for all the cases in Table 3.

We estimate the number of required function evaluations with multiple runs by multiplying the number for a single run by the multiple. By checking for a few cases we found that this estimate is accurate to within a single function evaluation. Comparing this table to the last two rows of Table 2 shows that multiple runs with a population of 2 and 2 DE iterations can cut down substantially the cost of reaching the optimum. With 6 multiple runs we can achieve a failure rate of 13.3% at a cost of 262 compared to 16% at a cost of 1013 without multiple runs. With 10 multiple runs costing 438 function evaluations we achieve a failure rate of about 4% that requires 2002 function evaluations without multiple runs.

*Table 3 Performance with multiple runs for DE/BFGS with a population of 2 and 2 DE iterations. Results are based on 100,800 replicate runs. Bump of 10/e applied to Optimum 1.*

| m | 1 | 2 | 3 | 4 | 5 | 6 | 7 | 8 | 9 | 10 |
|---|---|---|---|---|---|---|---|---|---|---|
| Percent failures | 71.3 | 51.0 | 36.2 | 26.0 | 18.4 | 13.3 | 9.3 | 6.9 | 4.7 | 3.6 |
| Percent expected if independent | 71.3 | 50.9 | 36.3 | 25.9 | 18.5 | 13.2 | 9.4 | 6.7 | 4.8 | 3.4 |
| Estimated number of function evaluations | 43.7 | 87.4 | 131 | 175 | 219 | 262 | 306 | 350 | 394 | 438 |

The table also shows that for this case the assumption of independent runs works well. This does not mean that it would work well for all problems, but it is encouraging to see this positive result. It turns out that while the independence assumption works well also for bumps in the other optima. Table 4 and Table 5 show the case of Table 3 when the bump is at Optimum 2, and Optimum 3, respectively. While adding bumps at these locations is slightly easier optimization problem than the bump at Optimum 1, it is seen that it is possible to achieve low probability of failure with multiple runs.

Table 4 Performance with multiple runs for DE/BFGS with a population of 2 and 2 DE iterations. Results are based on 100,800 replicate runs. Bump of 10/e applied to Optimum 2

| m | 1 | 2 | 3 | 4 | 5 | 6 | 7 | 8 | 9 | 10 |
|---|---|---|---|---|---|---|---|---|---|---|
| Percent failures | 62.0 | 38.3 | 23.9 | 14.8 | 9.0 | 5.8 | 3.5 | 2.2 | 1.3 | 0.74 |
| Percent expected if independent | 62.0 | 38.4 | 23.8 | 14.8 | 9.1 | 5.6 | 3.5 | 2.2 | 1.3 | 0.84 |
| Estimated number of function evaluations | 44.9 | 89.8 | 134.7 | 179.6 | 224.5 | 269 | 314 | 359 | 404 | 449 |

Table 5 Performance with multiple runs for DE/BFGS with a population of 2 and 2 DE iterations. Results are based on 100,800 replicate runs. Bump of 10/e applied to Optimum 3

| m | 1 | 2 | 3 | 4 | 5 | 6 | 7 | 8 | 9 | 10 |
|---|---|---|---|---|---|---|---|---|---|---|
| Percent failures | 65.4 | 42.6 | 27.6 | 18.2 | 11.7 | 7.5 | 5.1 | 3.4 | 2.1 | 1.4 |
| Percent expected if independent | 65.4 | 42.7 | 27.9 | 18.3 | 11.9 | 7.8 | 5.1 | 3.3 | 2.2 | 1.4 |
| Estimated number of function evaluations | 45 | 90 | 235 | 180 | 225 | 270 | 315 | 360 | 405 | 450 |

## VII. Concluding Remarks

Global optimization algorithms are often tested on easy functions with multiple identical-valued optima. This paper suggests that these functions could be fortified to become much harder by adding or subtracting a bump to one of the optima. This is illustrated for the Branin-Hoo function, which has three identical-valued global optima. Bumps were created using RBFs within a small vicinity of an optimum. The performance of the Python SciPy popular differential evolution (DE) optimizer, with and without a follow up by the BFGS local gradient based algorithm, was used on the fortified Branin-Hoo function. It is shown that adding the bump to the first optimum increases the required computational effort by about a factor of 16. Interestingly, it is found that higher chances of finding the optimum are afforded by repeating runs multiple times in comparison with running one optimization for more function evaluations. When a bump was added to Optimum 1, one optimization with 400 function evaluations resulted in about a 60% probability that the optimization failed to find the global optima. However, with nine optimizations each limited to 45 function evaluations (for a total of 405 function evaluations), the probability that the optimization failed was reduced to less than a 2%. For the first optimum, this reduced the computational effort to about a factor of 10 over the original Branin-Hoo (without the bump). Furthermore, it was found that the assumption that the multiple runs were independent gave reasonable estimates of the probability of finding the global optimum.